# ALL-PENTAGONAL FACE MULTI TORI

Mircea V. Diudea[a] and Aleksandar Ilić[b,*]

[a]Faculty of Chemistry and Chemical Engineering, "Babes-Bolyai" University,
Arany Janos Str. 11, 400028, Cluj, Romania
diudea@gmail.com

[b]Faculty of Sciences and Mathematics, University of Niš,
Višegradska 33, 18000 Niš, Serbia
aleksandari@gmail.com

**Abstract**. Design of monomeric unit, by using sequences of map operations, and the growing process in the building of a double-shell multi torus with all-pentagonal faces, is presented. It is shown that the monomer and some small intermediates, as hydrogenated species, have a moderate stability, between adamantane and $C_{60}$ fullerene, as calculated at the PM3 level of theory. The topology of these high genera structures is described in terms of Omega polynomial as a function of the net parameters. Close formulas for this polynomial and examples are presented.

## 1. Introduction

There are rigid monomers that can assembly in dendrimers, but the growth goes up only at the first generation. At the second generation, yet the endings of units are no more free, they fit to each other, thus forming either infinite lattice, if the unit symmetry is octahedral or spherical multi-torus, if the unit symmetry is tetrahedral.[1] The last one is the case of structures herein discussed.

Recall, dendrimers are hyper-branched nano-structures with rigorously tailored architecture. They can be functionalized at terminal branches, thus finding a broad pallet of applications in chemistry, medicine, etc. [2-5] Molecular topology of dendrimers is a continuously developed topic. It includes vertex and fragment enumeration as well as calculation of some topological descriptors, such as topological indices, sequences of numbers or polynomials.[5-8]

Multi-tori are complex structures consisting of more than one tubular ring. They are supposed to result by self-assembly of some repeat units or monomers which can be designed by opening of cages/fullerenes or by appropriate map/net operations. Such structures can appear in spongy carbon or in zeolites.[9,10]

Structures of high genera can be designed starting from the Platonic solids, by using appropriate map operations. Such units may form either infinite lattices of negative curvature or closed cages, showing high porosity. Such structures have also been modeled by Lenosky et al.,[11] Terrones and Mackay,[12,13] etc. Spongy carbons have been recently synthesized.[14-16]

## 2. Cage Building

The tetrapodal monomer *tt* designed by $Trs(P_4(T))$ sequence of map operations[1] forms a dendrimer (Figure 1), which at the second generation, of 17 units, shows pentagonal cycles/ super-faces appear (Figure 2, left). The "growing process" was leaded so that minimum number of monomers *tt* is added to close a super-cycle. At $M_{27}$, a spherical unit *u* already appears (Figure 2, right). The unit $U_1$ (Figure 3, left) is in fact a multi torus $MT(Trs(P_4(T)))\&[S_2(D)]$, with the core designed by the $S_2$ operation applied to the Dodecahedron D. It can evolve either linearly (Figure 3, right and Figure 4, right) or arrange as cyclic (Figure 5, left) or multi-shell multi tori, like MT(12U) (Figure 5, right).

These are structures of high genera: $U_1$ is $g=21$, the highest genus is MT(12U): $g=131$. These numbers resulted by applying the Euler's formula:[17] $v - e + f = 2(1 - g)$, where $v$, $e$, $f$ and $g$ are the number of vertices/atoms, edges/bonds, faces and genus, respectively. The genus is the number of handles required to be attached to the sphere to make it homeomorphic to a given embedding of a graph.[18] It is not trivial to count the number of simple tori (*i.e.*, the crude meaning of the genus $g$) in such complex structures.

The number of monomers/*tt* in MT(12U), all-pentagonal face multi-torus, is 130, that fits the number of atoms/points in the super-array of 12 Dodecahedra[1,9,10,16] 12D; it means every point was substituted by a tetrapodal unit *tt*. On the other hand, the number of hollows $R_{10}$ is exactly 114, as the number of faces in 12D.

$M_1$; $v=22$; $e=36$; $f_5=12$; $tt=1$; $g=2$     $M_5$; $v=98$; $tt=5$

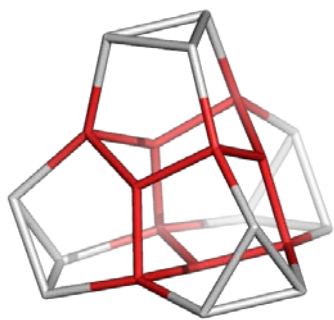 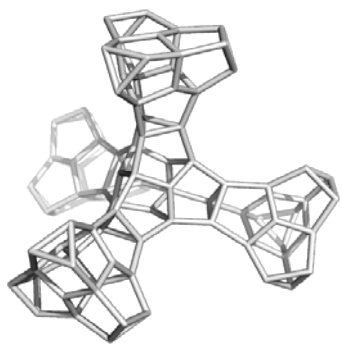

Figure 1. Monomer designed by $Trs(P_4(T))$ and a first generation dendrimer

$M_{17}\_6r$; $v=308$; $tt=17$            $M_{27}\_1u$; $v=474$; $tt=27$

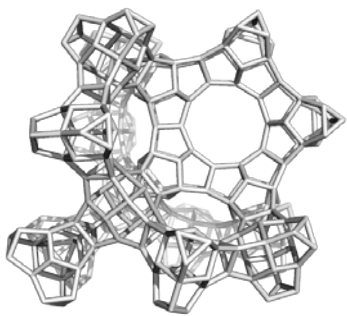 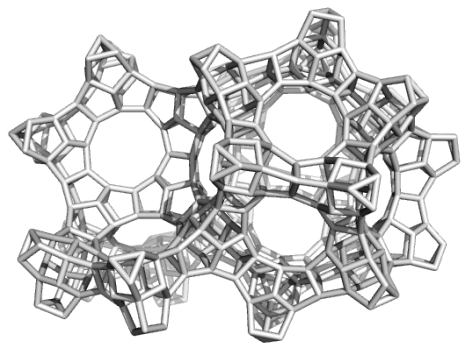

Figure 2. Growing process: apparition of cyclic and spherical units

$U_1$; $v=350$; $e=630$; $f_5=240$; $tt=20$; $g=21$; five-fold symmetry        $U_2$; $v=605$; $tt=35$

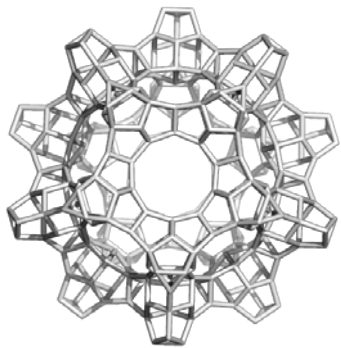 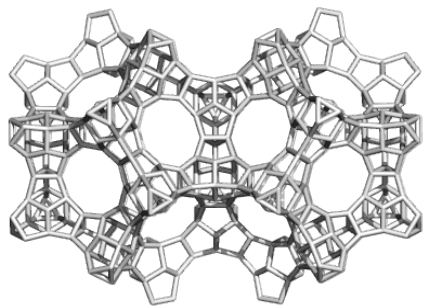

Figure 3. Spherical unit $U_1 = MT(Trs(P_4(T)))\&[S_2(D)]$ (left) and its dimer (right)

$M_{34}\_1u$; $v=598$; $tt=34$            $U_3$; $v=860$; $tt=50$

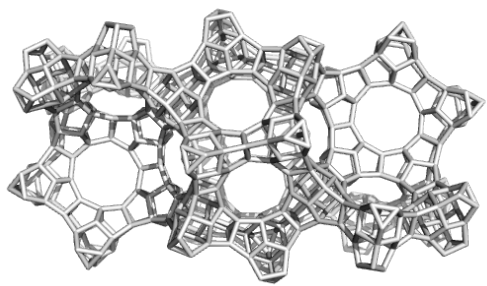 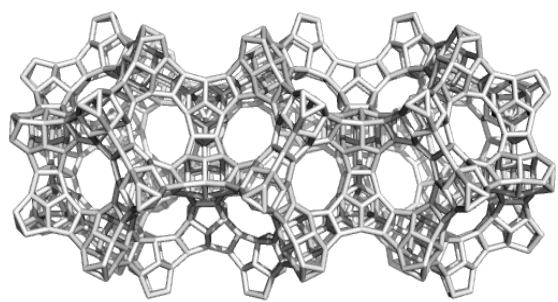

Figure 4. Growing process: apparition of spherical units (left) and a linear array of the spherical units.

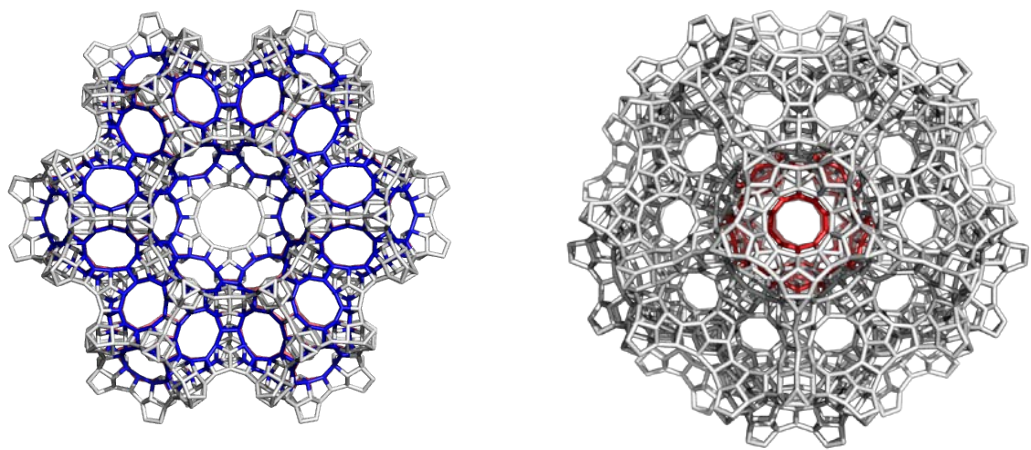

$U_{6,cyc}$; $v=1530$; $tt=90$     MT(12U); $v=2170$; $e=3990$; $f_5=1560$; $R_{10}=114$; $tt=130$; $g=131$; five-fold symmetry

Figure 5. Cyclic (left) and multi-shell spherical (right) arrays of $U_1$.

### 3. Stability of Tetrapodal Structures

Hydrogenated structures involved in these molecular constructions have shown moderate stability as given by their heat of formation HF, total energy TE and HOMO-LUMO Gap HLGAP, calculated at the PM3 level of theory (Table 1). For example, the total energy per heavy atoms of the structures in Table 1 are between the values of adamantane (-3305.19 kcal/mol), which is the most related small structure (see Figure 1, left, in red) and $C_{60}$ (-2722.45 kcal/mol), the standard molecule in nanostructures. The same is true about the HOMO-LUMO gap. Calculations by using a density functional-based tight binding method combined with the self-consistent charge technique (SCC-DFTB) on hydrogenated units of diamond and a diamond-like network[19] have shown the same ordering of stability as given by PM3 approach; thus, our results reported here can be considered as pertinent ones.

Table 1. Quantum Chemistry data for some units designed by $Trs(P_4(T))$:
Heat of Formation HF, Total energy TE and HOMO-LUMO Gap HLGAP

| Name | N-heavy atoms | HF (kcal/mol) | HF/N heavy | TE (kcal/mol) | TE/N heavy | HLGAP (eV) | Sym. |
|---|---|---|---|---|---|---|---|
| $M_1$ | 22 | 308.484 | 14.022 | -65540 | -2979.09 | 11.9925 | $t_d$ |
| $M_5$ | 98 | 1510.296 | 15.411 | -286415 | -2922.6 | 11.6847 | $c_1$ |
| $M_{17}$ | 308 | 5095.126 | 16.543 | -887083 | -2880.14 | 11.5121 | $c_1$ |

## 4. Omega Polynomial

In a connected graph $G(V,E)$, with the vertex set $V(G)$ and edge set $E(G)$, two edges $e = uv$ and $f = xy$ of $G$ are called *codistant e co f* if they obey the relation:[20]

$$d(v,x) = d(v,y) + 1 = d(u,x) + 1 = d(u,y) \quad (1)$$

which is reflexive, that is, *e co e* holds for any edge *e* of $G$, and symmetric, if *e co f* then *f co e*. In general, relation *co* is not transitive; if "*co*" is also transitive, thus it is an equivalence relation, then $G$ is called a *co-graph* and the set of edges $C(e) := \{f \in E(G); f\ co\ e\}$ is called an *orthogonal cut oc* of $G$, $E(G)$ being the union of disjoint orthogonal cuts: $E(G) = C_1 \cup C_2 \cup ... \cup C_k$, $C_i \cap C_j = \emptyset, i \neq j$. Klavžar[21] has shown that relation *co* is a theta Djoković-Winkler relation.[22,23]

We say that edges *e* and *f* of a plane graph $G$ are in relation *opposite, e op f,* if they are opposite edges of an inner face of $G$. Note that the relation *co* is defined in the whole graph while *op* is defined only in faces. Using the relation *op* we can partition the edge set of $G$ into *opposite* edge *strips, ops*. An *ops* is a quasi-orthogonal cut *qoc*, since *ops* is not transitive.

Let $G$ be a connected graph and $s_1, s_2, ..., s_k$ be the *ops* strips of $G$. Then the *ops* strips form a partition of $E(G)$. The length of *ops* is taken as maximum. It depends on the size of the maximum fold face/ring $F_{max}/R_{max}$ considered, so that any result on Omega polynomial will have this specification.

Denote by $m(G,s)$ the number of *ops* of length $s$ and define the Omega polynomial as:[24-26]

$$\Omega(G,x) = \sum_s m(G,s) \cdot x^s \quad (2)$$

Its first derivative (in $x=1$) equals the number of edges in the graph:

$$\Omega'(G,1) = \sum_s m(G,s) \cdot s = e = |E(G)| \quad (3)$$

On Omega polynomial, the Cluj-Ilmenau[20] index, $CI=CI(G)$, was defined:

$$CI(G) = \{[\Omega'(G,1)]^2 - [\Omega'(G,1) + \Omega''(G,1)]\} \quad (4)$$

## 5. Omega Polynomial and Topology of MT(12U)

The simplicity and utility of Omega polynomial description, in these tetrapodal multi tori entirely tessellated by pentagons, is fully demonstrated.[1] At maximum face/ring $F_5$ calculation, the odd faces give a single term at exponent 1, the coefficient being the number of edges in the graph.

At maximum ring $R_6$, there are only two terms; the coefficient of term at the exponent 3 is four times the number of tetrapodes *tt*. The first derivative of this term equals the number of pentagonal faces in structure, which is also given by $12 \times tt$.

It is important for a description to be both diagnostic and prognostic, which the Omega polynomial is. Formulas (function of the net parameters) and examples for Omega polynomial and CI index (in parenthesis) are given in Table 2 and Table 3.

Numerical calculations have been done by our software package Nano Studio.[27]

Table 2. Omega Polynomial and Topology of Multi Tori

|  | Formulas |
|---|---|
| $F_5$ | $\Omega(x; F_5; M_{m,r}) = [36m - 3(m+r-1)] \cdot x^1; CI = 3(11m - r + 1)(33m - 3r + 2)$ |
|  | $\Omega(x; F_5; U_{u\text{-Lin}}) = [630u - 165(u-1)] \cdot x^1 = (465u + 165) \cdot x^1; CI = 15(31u + 11)(465u + 164)$ |
|  | $\Omega(x; F_5; U_{u\text{-Cyc}}) = (630u - 165u) \cdot x^1 = 465u \cdot x^1; CI = 465u(465u - 1)$ |
|  | $\Omega(x; F_5; MT(12U)) = 133d \cdot x^1; d = 30$ |
| Examples | $M_{1,0}$: $36x^1$ (1260); $M_{3,0}$: $102x^1$ (10302); $M_{4,0}$: $135 x^1$ (18090); $M_{5,0}$: $168x^1$ (28056) |
|  | $U_{1Lin}$: $630 x^1$ (396270); $U_{3Lin}$: $1560x^1$ (2432040); $U_{4Lin}$: $2025x^1$ (4098600) |
|  | $U_{6Cyc}$: $2790x^1$ (7781310) |
|  | MT(12U): $3990x^1$ (15916110) |
| $R_6$ | $\Omega(x; R_6; M_{m,r}) = [24m - 3(m+r-1)] \cdot x^1 + 4m \cdot x^3; CI = (33m - 3r + 3)^2 - 3(19m - r + 1)$ |
|  | $\Omega(x; R_6; U_{u\text{-Lin}}) = [390u - 105(u-1)] \cdot x^1 + [80u - 20(u-1)]x^3 =$ |
|  | $(285u + 105) \cdot x^1 + (60u + 20) \cdot x^3; CI = (465u + 165)^2 - 825u - 285$ |
|  | $\Omega(x; R_6; U_{u\text{-Cyc}}) = (390u - 105u) \cdot x^1 + (80u - 20u)x^3 = 285u \cdot x^1 + 60u \cdot x^3;$ |
|  | $CI = 75u(2883u - 11)$ |
|  | $\Omega(x; R_6; MT(12U)) = 81d \cdot x^1 + 26t \cdot x^3; d = 30; t = 20$ |
|  | *d* = number of U-double joints; *t* = number of U-triple joints. |
| Examples | $M_{1,0}$: $24x^1 + 4x^3$ (1236); $M_{5,0}$: $108x^1 + 20x^3$ (27936); $M_{17,6}$: $342x^1 + 68x^3$ (297162) |
|  | $U_{1Lin}$: $390x^1 + 80x^3$ (395790); $U_{3Lin}$: $960x^1 + 200x^3$ (2430840) $U_{4Lin}$: $1245x^1 + 260x^3$ (4097040) |
|  | $U_{6Cyc}$: $1710x^1 + 360x^3$ (7779150) |
|  | MT(12U): $2430x^1 + 520x^3$ (15912990) |

Table 3. The number of tetrapodes and the number of atoms of Multi Tori

| | |
|---|---|
| $tt$ | $tt = 20u - (5d - 2t)$ |
| | $tt_{\text{U-Lin}} = 20u - 5(u-1)$ |
| | $tt_{\text{U-Cyc}} = 20u - 5u$ |
| Examples | $tt_{3-\text{Lin}} = 50;\ tt_{4-\text{Lin}} = 65$ |
| | $tt_{6-\text{Cyc}} = 90$ |
| | $tt_{\text{MT}} = 20 \times 12 - (5 \times 30 - 2 \times 20) = 130$ |
| $v$ | $v(\text{M}_m) = 19m + 3;\ m = 1, 2, .., 11;\ v(\text{M}_m) = 16m + 36;\ m = 12, .., 17$ |
| | $v(\text{U}_{u-\text{Lin}}) = 255u + 95;\ u = 1, 2, ...$ |
| | $v(\text{U}_{u-\text{Cyc}}) = 256u;\ u = 6, 7, ...$ |
| | $v(\text{MT}(12\text{U})) = 5(5tt - 18u);\ tt = 130;\ u = 12$ |
| Examples | $\text{M}_1$: 22; $\text{M}_3$: 60; $\text{M}_4$: 70; $\text{M}_5$: 98; $\text{M}_{12}$: 228; $\text{M}_{15}$: 276; $\text{M}_{17}$: 308 |
| | $\text{U}_1$: 350; $\text{U}_3$: 860; $\text{U}_4$: 1115 |
| | $\text{U}_{6-\text{Cyc}}$: 1536 |
| | MT(12U): 2170 |

## 5. Conclusions

Design of a hypothetical dendrimer and multi-tori by using the sequence of map operation $Op(Trs(P_4(\text{T})))$ was presented. The computed total energy of the involved repeat unit and some small intermediates showed stability intermediate to that of adamantane and $C_{60}$, the reference molecule in nanostructures. The topology of the proposed structures was described in terms of Omega polynomial as a function of the net parameters.

**Acknowledgements**: The first author (MVD) is supported by the Romanian CNCSIS-UEFISCSU project number PN-II IDEI 129/2010, and the second author (AI) is supported by the Research Grant 144007 of Serbian Ministry of Science.

**References**

1. M. V. Diudea, *Nanomolecules and Nanostructures - Polynomials and Indices*, Mathematical Chemistry Monographs 10, University of Kragujevac, Serbia, 2010.
2. D. A. Tomalia, A. M. Naylor, W. A. I. Goddard, Starburst Dendrimers: Molecular- Level Control of Size, Shape, Surface Chemistry, Topology, and Flexibility from Atoms to Macroscopic Matter, *Angew. Chem. Int. Ed.* 1990, *29*, 138-175.
3. G. R. Newkombe, C. N. Moorefield, F. Voegtle, *Dendrimers and dendrons: Concepts, syntheses, applications*; Wiley-VCH, Weinheim, 2001.

4. D. A. Tomalia, J. M. J. Frechet, Discovery of dendrimers and dendritic polymers: A brief historical perspective, *J. Polym. Sci.*, *Part A: Polym. Chem.* 2002, *40*, 2719-2728

5. M. V. Diudea, G. Katona, in: G. A. Newkome, Ed., *Advan. Dendritic Macromol.*, 1999, *4*, 135-201.

6. M. V. Diudea, A. A. Kiss, E. Estrada, N. Guevara, Connectivity-, Wiener- and Harary-type indices of dendrimers, *Croat. Chem. Acta*, 2000, *73*, 367-381.

7. M. V. Diudea, Wiener Index of Dendrimers, *MATCH Commun. Math. Comput. Chem.*, 1995, *32*, 71-83.

8. M. V. Diudea, B. Parv, Hyper-Wiener Index of Dendrimers, *J. Chem. Inf. Comput. Sci.* 1995, *35*, 1015-1018.

9. M. V. Diudea, (Ed.), *Nanostructures, Novel Architecture*, NOVA, New York, 2005.

10. M. V. Diudea, Cs. L. Nagy, *Periodic Nanostructures*, Springer, 2007.

11. T. Lenosky, X. Gonze, M. Teter, V. Elser, Energetics of negatively curved graphitic carbon, *Nature*, 1992, *355*, 333-335.

12. H. Terrones, A. L. Mackay, From $C_{60}$ to negatively curved graphite, *Prog. Crystal Growth and Charact.*, 1997, *34*, 25-36.

13. H. Terrones, A. L. Mackay, Triply periodic minimal surfaces decorated with curved graphite, *Chem. Phys. Lett.*, 1993, *207*, 45-50.

14. E. Barborini, P. Piseri, P. Milani, G. Benedek, C. Ducati, J. Robertson, Negatively curved spongy carbon, *Appl. Physics Lett.*, 2002, *81*, 3359-3361.

15. G. Benedek, H. Vahedi-Tafreshi, E. Barborini, P. Piseri, P. Milani, C. Ducati, J. Robertson, The structure of negatively curved spongy carbon, *Diamond and Related Materials*, 2003, *12*, 768–773.

16. Cs. L. Nagy, M. V. Diudea, Nanoporous Carbon Structures, in: M. V. Diudea, Ed., *Nanostructures-Novel Architecture*. NOVA, New York, 2005, 311-334.

17. L. Euler, Elementa doctrinae solidorum, *Novi Commentarii Academiae Scientiarum Imperialis Petropolitanae*, 1758, vol. 4, 109-160.

18. F. Harary, *Graph Theory*, Addison-Wesley, Reading, MA, 1969.

19. M. V. Diudea, A. Bende, D. Janežić, Omega polynomial in diamond-like networks, *Fullerenes, Nanotubes, Carbon Nanostruct.*, 2010, *18* (3), 000-000.

20. P. E. John, A. E. Vizitiu, S. Cigher, M. V. Diudea, CI Index in Tubular Nanostructures, *MATCH Commun. Math. Comput. Chem.*, 2007, *57*, 479-484.

21. S. Klavžar, Some comments on co graphs and CI index, *MATCH Commun. Math. Comput. Chem.*, 2008, *59*, 217-222.

22. D. Ž. Djoković, Distance preserving subgraphs of hypercubes, *J. Combin. Theory Ser. B*, 1973, *14*, 263-267.

23. P. M. Winkler, Isometric embedding in products of complete graphs, *Discrete Appl. Math.*, 1984, *8*, 209-212.

24. M. V. Diudea, Omega Polynomial, *Carpath. J. Math.*, 2006, *22*, 43-47.

25. M. V. Diudea, S. Cigher, P. E. John, Omega and Related Counting Polynomials, *MATCH Commun. Math. Comput. Chem.*, 2008, *60*, 237-250.

26. M. V. Diudea, S. Cigher, A. E. Vizitiu, M. S. Florescu, P. E. John, *J. Math. Chem.*, 2009, *45*, 316-329.

27. Cs. L. Nagy, M. V. Diudea, *Nano Studio software*, Babes-Bolyai University, 2009.